\documentclass{amsart}

\usepackage{amsmath,amsthm,amsfonts,amscd,amssymb}
\usepackage{verbatim}
\usepackage[pdfborder={0 0 0}]{hyperref}
\usepackage{color} 

\newtheorem{thm}{Theorem}[section]
\newtheorem{cor}[thm]{Corollary}
\newtheorem{lem}[thm]{Lemma}

\theoremstyle{definition}

\theoremstyle{remark}
\newtheorem{rem}{Remark}[section]

\begin{document}

\title{Small representations of integers by integral quadratic forms}

\author{Wai Kiu Chan and Lenny Fukshansky}\thanks{The second author was partially supported by the Simons Foundation grant \#519058.}

\address{Department of Mathematics and Computer Science, Wesleyan University, Middletown, CT 06459}
\email{wkchan@wesleyan.edu}

\address{Department of Mathematics, 850 Columbia Avenue, Claremont McKenna College, Claremont, CA 91711}
\email{lenny@cmc.edu}


\begin{abstract} Given an isotropic quadratic form over a number field which assumes a value $t$, we investigate the distribution of points at which this value is assumed. Building on the previous work about the distribution of small-height zeros of quadratic forms, we produce bounds on height of points outside of some algebraic sets in a quadratic space at which the form assumes the value $t$.   Our bounds on height are explicit in terms of the heights of the form, the space, the algebraic set and the value $t$.
\end{abstract}

\maketitle

\def\A{{\mathcal A}}
\def\AA{{\mathfrak A}}
\def\B{{\mathcal B}}
\def\C{{\mathcal C}}
\def\D{{\mathcal D}}
\def\E{{\mathcal E}}
\def\F{{\mathcal F}}
\def\Ff{{\mathfrak F}}
\def\G{{\mathcal G}}
\def\x{{\mathcal H}}
\def\I{{\mathcal I}}
\def\J{{\mathcal J}}
\def\K{{\mathcal K}}
\def\kk{{\mathfrak K}}
\def\L{{\mathcal L}}
\def\LL{{\mathfrak L}}
\def\M{{\mathcal M}}
\def\O{{\mathcal O}}
\def\W{{\omega}}
\def\CC{{\mathfrak C}}
\def\mm{{\mathfrak m}}
\def\MM{{\mathfrak M}}
\def\OO{{\mathfrak O}}
\def\P{{\mathcal P}}
\def\R{{\mathcal R}}
\def\s{{\mathcal S}}
\def\V{{\mathcal V}}
\def\X{{\mathcal X}}
\def\XX{{\mathfrak X}}
\def\Y{{\mathcal Y}}
\def\Z{{\mathcal Z}}
\def\H{{\mathcal H}}
\def\cee{{\mathbb C}}
\def\pee{{\mathbb P}}
\def\que{{\mathbb Q}}
\def\real{{\mathbb R}}
\def\zed{{\mathbb Z}}
\def\hyp{{\mathbb H}}
\def\aaa{{\mathbb A}}
\def\ff{{\mathbb F}}
\def\kk{{\mathfrak K}}
\def\qbar{{\overline{\mathbb Q}}}
\def\kbar{{\overline{K}}}
\def\ybar{{\overline{Y}}}
\def\kkbar{{\overline{\mathfrak K}}}
\def\ubar{{\overline{U}}}
\def\eps{{\varepsilon}}
\def\ahat{{\hat \alpha}}
\def\bhat{{\hat \beta}}
\def\gt{{\tilde \gamma}}
\def\h{{\tfrac12}}
\def\dd{{\partial}}
\def\baa{{\boldsymbol \alpha}}
\def\bfa{{\boldsymbol a}}
\def\bfb{{\boldsymbol b}}
\def\be{{\boldsymbol e}}
\def\bei{{\boldsymbol e_i}}
\def\bff{{\boldsymbol f}}
\def\bc{{\boldsymbol c}}
\def\bm{{\boldsymbol m}}
\def\bk{{\boldsymbol k}}
\def\bi{{\boldsymbol i}}
\def\bl{{\boldsymbol l}}
\def\bq{{\boldsymbol q}}
\def\bu{{\boldsymbol u}}
\def\bt{{\boldsymbol t}}
\def\bs{{\boldsymbol s}}
\def\bfu{{\boldsymbol u}}
\def\bv{{\boldsymbol v}}
\def\bw{{\boldsymbol w}}
\def\bx{{\boldsymbol x}}
\def\bX{{\boldsymbol X}}
\def\bz{{\boldsymbol z}}
\def\bwy{{\boldsymbol y}}
\def\bY{{\boldsymbol Y}}
\def\bL{{\boldsymbol L}}
\def\ba{{\boldsymbol a}}
\def\bb{{\boldsymbol b}}
\def\bet{{\boldsymbol\eta}}
\def\bxi{{\boldsymbol\xi}}
\def\bo{{\boldsymbol 0}}
\def\bol{{\boldkey 1}_L}
\def\ep{\varepsilon}
\def\p{\boldsymbol\varphi}
\def\q{\boldsymbol\psi}
\def\rank{\operatorname{rank}}
\def\aut{\operatorname{Aut}}
\def\lcm{\operatorname{lcm}}
\def\sgn{\operatorname{sgn}}
\def\spn{\operatorname{span}}
\def\md{\operatorname{mod}}
\def\Norm{\operatorname{Norm}}
\def\dim{\operatorname{dim}}
\def\det{\operatorname{det}}
\def\Vol{\operatorname{Vol}}
\def\rk{\operatorname{rk}}
\def\ord{\operatorname{ord}}
\def\ker{\operatorname{ker}}
\def\div{\operatorname{div}}
\def\Gal{\operatorname{Gal}}
\def\GL{\operatorname{GL}}
\def\p{\operatorname{p}}
\def\q{\operatorname{q}}
\def\t{\operatorname{t}}
\def\hs{{\hat \sigma}}
\def\chr{\operatorname{char}}

\section{Introduction and statement of results}

The question of effective distribution of zeros of quadratic forms has a long history. In particular, the famous Cassels' theorem (1955) asserts that an isotropic rational quadratic form has zeros of bounded height with effective bounds depending on the height of this form,~\cite{cassels:small}. By height here we mean a height function, a standard tool of Diophantine geometry, measuring the arithmetic complexity of an object; we provide a brief overview of height functions in Section~\ref{tools}. There have been many extensions and generalizations of Cassels theorem over the years, see~\cite{cassels_overview} for a survey of this lively area. In particular, in the recent years there have been several results assessing the distribution of small-height zeros by proving that they cannot be easily ``cut-out" by any finite union of linear subspaces or even projective varieties not completely containing the quadratic hypersurface in question. Specifically, results of \cite{smallzeros}, \cite{dietmann}, \cite{cfh}, \cite{GR} prove existence of zeros of bounded height of an isotropic quadratic form outside of a union of projective varieties. One may wonder if the same type of results hold for other values of a quadratic form. Namely, if $K$ is a number field, $Q$ is a nonzero quadratic form on a subspace $V$ of $K^n$ and $t$ is an element of $K$ represented by $(V, Q)$, can the small-height points $\bz \in V$ with $Q(\bz)= t$ be efficiently ``cut-out" by sufficiently simple algebraic sets? This question can be viewed as an inhomogeneous analogue of the problem about small-height zeros. To state our results, let us start with some notation.

Let $K$ be a number field and $\O_K$ its ring of integers. Let $j \geq 1$ be an integer. For each $1 \leq i \leq j$, let $\s_i$ be a finite set of polynomials in $K[X_1, \ldots, X_n]$ and $Z_K(\s_i)$ be its zero set in $K^n$, that is,
$$Z_K(\s_i) = \{\bx \in K^n : P(\bx) = 0 \mbox{ for all } P \in \s_i\}.$$
For the collection $\bs: = \{ \s_1, \ldots, \s_j \}$ of finite sets of polynomials, define
\begin{equation}
\label{Z_K}
\Z_{\bs} := \bigcup_{i=1}^j Z_K(\s_i),
\end{equation}
and
\begin{equation}
\label{def_M}
M_{\bs} := \sum_{i=1}^j \max \{\deg P: P \in \s_i\}.
\end{equation}
The set $\Z_\bs$ defined in \eqref{Z_K} is said to be {\em homogeneous} if all the polynomials in the sets $\s_1, \ldots, \s_j$ are homogeneous polynomials.  When $\Z_\bs$ is not homogeneous, we will sometimes need the {\it homogenization} of $\Z_{\bs}$, defined as follows. Introducing one more variable $X_{n+1}$, we define the homogenization $P^*(X_1,.\dots,X_{n+1}) \in K[X_1,\dots,X_{n+1}]$ of each polynomial $P(X_1,\dots,X_n) \in K[X_1,\dots,X_n]$ to be the unique homogeneous polynomial such that $\deg P^* = \deg P$ and
$$P(X_1,\dots,X_n) = P^*(X_1,\dots,X_n,1),$$
and for each $1 \leq i \leq j$, let $\s^*_i = \{ P^* : P \in \s_i \}$. Then define $\bs^* = \{ \s^*_1, \ldots, \s^*_j \}$ and, as above, let  $Z_K(\s^*_i) = \{\bx \in K^{n+1} : P^*(\bx) = 0 \mbox{ for all } P \in \s_i\}$, as well as
\begin{equation}
\label{Z*_K}
\Z_{\bs^*} := \bigcup_{i=1}^j Z_K(\s^*_i).
\end{equation}
Then $\Z_{\bs^*}$ is a homogeneous algebraic set in $K^{n+1}$, while clearly
$$M_{\bs^*} := \sum_{i=1}^j \max \{\deg P^* : P \in \s_i\} = M_{\bs}.$$
We also write $H$ and $h$ for the homogeneous and inhomogeneous height function, respectively, to be defined in Section~\ref{tools}. Further, if $V \subseteq K^n$ and $Q(X_1, \ldots, X_n) \in K[X_1, \ldots, X_n]$ is a quadratic form in $n$ variables, we refer to the pair $(V,Q)$ as a
quadratic space (over $K$).  This quadratic space is called isotropic if there exists $\bo\neq \bx \in V$ such that $Q(\bx) = 0$.  A totally isotropic subspace of $V$ is a subspace of $V$ on which $Q$ vanishes identically.  A maximal totally isotropic subspace of $V$ is the largest totally isotropic subspace of $V$ with respect to inclusion.   The orthogonal complement of $V$ in $K^n$ is defined to be
$$\perp_Q(V) : = \{\bx \in K^n : Q(\bx, \bwy) = 0\ \mbox{ for all } \bwy \in V \},$$
where $Q(\bx, \bwy)$ is the associated symmetric bilinear form.  It is easy to see that $\perp_Q(V)$ is a subspace of $K^n$. The radical $V^\perp$ of $(V,Q)$ is defined to be the subspace of singular points in $V$ with respect to $Q$, i.e.
$$V^\perp := \{\bx \in V : Q(\bx, \bwy) = 0 \mbox{ for all } \bwy \in V \}.$$
Then $\dim(V) + \dim(\perp_Q(V)) - \dim(V^\perp) = n$.  The quadratic space $(V,Q)$ is called regular if $V^\perp = \{\bo\}$; we may also say that $Q$ is regular on $V$, or simply that $Q$ is regular if $V = K^n$. The dimension of a maximal totally isotropic subspace of a regular quadratic space $V$ is called the Witt index of $V$.  With this notation and terminology, here is a simplified version of a result of~\cite{GR} (which is an improvement of the previous similar results \cite{cfh}, \cite{dietmann}, \cite{smallzeros}), adapted to our purposes.

\begin{thm} [\cite{GR}, Theorem~7.1] \label{miss_hyper} Let $Q$ be a nonzero quadratic form in $n$ variables over $K$, $V$ be an $m$-dimensional subspace of $K^n$, and $w$ be the dimension of a maximal totally isotropic subspace of the quadratic space $(V,Q)$. Let $\Z_{\bs}$ be a homogeneous algebraic set as in~\eqref{Z_K} such that $Q$ has a nontrivial zero in $V \setminus \Z_{\bs}$.  Then there exists a nontrivial zero $\bx \in \left( V \cap \O_K^n \right) \setminus \Z_{\bs}$ of $Q$ such that
\begin{equation}
\label{miss_hyper_bnd}
H(\bx) \ll H(Q)^{\frac{m-w+1}{2}} H(V)^2,
\end{equation}
where the implied constant depends only on $K$, $m$ and $M_{\bs}$.
\end{thm}

\begin{rem} Here and in the rest of the paper we choose not to explicitly specify the dimensional and field constants: their importance is marginal while their form often gets complicated. This being said, all the results of this paper are effective in that these constants can be explicitly computed.
\end{rem}

Results like Theorem~\ref{miss_hyper} are usually stated with the point $\bx$ in question having coordinates in $K$. It is however always possible to find some $a \in \O_K$ such that $a\bx \in \O_K^n$; since $Q(a\bx) = a^2 Q(\bx) = 0$ and the homogeneous height $H(a\bx) = H(\bx)$, we can simply assume that $\bx \in \O_K^n$. On the other hand, the distinction between points in $K^n$ and $\O_K^n$ is more subtle in the inhomogeneous situations.  For any $0 \neq t \in Q(V)$, let $Q_t(X_1, \ldots, X_n): = Q(X_1, \ldots, X_n) - t$ and $Q_t^*(X_1, \ldots, X_{n+1})$ be its homogenization as defined earlier.    Here is the first ``inhomogeneous" observation, which we derive from Theorem~\ref{miss_hyper} in Section~\ref{field}.

\begin{thm} \label{QZ}  Let $Q$ be a nonzero quadratic form in $n$ variables over $K$ and let $V \subseteq K^n$ be an $m$-dimensional subspace, $m \leq n$.  Let $w \geq 0$ be the dimension of a maximal totally isotropic subspace of the quadratic space $(V,Q)$. Let $0 \neq t \in Q(V)$ and $\Z_{\bs}$ be an algebraic set as in~\eqref{Z_K} which does not contain the zero set of the polynomial $Q_t(X_1, \ldots, X_n)$. Then there exists a point $\bz \in V \setminus \Z_{\bs}$ such that $Q(\bz) = t$ and
\begin{equation}
\label{QZ-1}
h(\bz) \ll H(Q_t^*)^{\frac{m-w+2}{2}} H(V)^2.
\end{equation}
The implied constant in the inequality depends only  on $K$, $m$, $n$ and $M_{\bs}$.
\end{thm}

\noindent
The idea of the proof here is to consider $Q_t^*(X_1, \ldots, X_{n+1})$ under the condition $X_{n+1} \neq 1$, apply Theorem~\ref{miss_hyper} and then divide out by this new variable.

One conceptual difference between Theorem~\ref{miss_hyper} and Theorem~\ref{QZ} is that in the homogeneous case the small-height point in question has coordinates in~$\O_K$ whereas in the inhomogeneous setting it does not. In fact, the problem of obtaining an inhomogeneous result over $\O_K$ instead of~$K$ is substantially more difficult. Some state-of-the art results for integer solutions to inhomogeneous integral quadratic equations without any avoidance conditions can be found in \cite{dietmann2}. Theorem \ref{smallrep} below is a version of the main result of \cite{dietmann2}, which addresses the special case $K = \mathbb Q$, stated in terms of our notation.  A quadratic form in $\mathbb Z[X_1, \ldots, X_n]$ is said to be integral if its Gram matrix has integer entries.

\begin{thm} [\cite{dietmann2}, Theorem 1] \label{smallrep}
Let $Q(X_1, \ldots, X_n) \in \mathbb Z[X_1, \ldots, X_n]$ be a regular integral quadratic form.  Suppose that an integer $t$ is integrally represented by $Q$.  Then there is a vector $\bz \in \mathbb Z^n$ such that $Q(\bz) = t$ and
$$h(\bz) \ll h(Q_t)^{\ell(n)},$$
where the implied constant depends only on $n$ and
$$\ell(n) = \begin{cases}
2100 & \mbox{ if $n = 3$},\\
84 & \mbox{ if $n = 4$},\\
5n + 19 + 74/(n - 4) & \mbox{ if $n \geq 5$}.
\end{cases}$$
\end{thm}

The following corollary is a subspace version of Theorem \ref{smallrep}.  We postpone its proof to the next section after  all the height functions are defined.

\begin{cor} \label{smallrepV}
Let $Q(X_1, \ldots, X_n) \in \mathbb Z[X_1, \ldots, X_n]$ be an integral quadratic form and $(V, Q)$ be an $m$-dimensional regular subspace of $\mathbb Q^n$, $m \geq 3$. Suppose that $t \in Q(V\cap \mathbb Z^n)$.  Then there exists $\bz \in V \cap \mathbb Z^n$ such that $Q(\bz) = t$ and
\begin{equation} \label{smallrepVbound}
h(\bz) \ll h(Q_t)^{\ell(m)} H(V)^{2\ell(m) + 1}.
\end{equation}
\end{cor}

As for integer solutions to inhomogeneous integral quadratic equations with avoidance conditions, we have the following result in this direction.  Its proof will be given in Section \ref{proofs}.

\begin{thm} \label{subspaces} Let $Q(X_1, \ldots, X_n)\in \mathbb Z[X_1, \ldots, X_n]$ be an integral quadratic form. Let $V \subseteq \que^n$ be an $m$-dimensional subspace, $3 \leq m \leq n$, such that the quadratic space $(V,Q)$ is regular and isotropic of Witt index $w$.  Let $W_1, \ldots, W_k$ be distinct hyperplanes of $V$.  Then for any $0 \neq t \in Q(V\cap \zed^n)$,  there exists $\bz \in \left( V \cap \zed^n \right) \setminus \bigcup_{i=1}^k W_i$ such that $Q(\bz) = t$ and
$$h(\bz) \ll  \begin{cases}
h(Q_t)^{(1 + \frac{2}{m-2})\ell(m) + m + 2 + \frac{m+4}{m-2}} H(V)^{(2 + \frac{4}{m-2})\ell(m) + 5 + \frac{8}{m-2}} & \mbox{ if $w = 1$};\\
h(Q_t)^{2\ell(m) + \frac{m-w+5}{2}}H(V)^{2\ell(m) + 3} & \mbox{ if $w \geq 2$},
\end{cases} $$
where the implied constant depends only on $m$, $n$, and $k$.
\end{thm}

\begin{rem} Notice that Theorem~\ref{subspaces} cannot be true when the dimension of $V$ is equal to~$2$, since in that case the number of representations of an integer by $Q$ on $V \cap \zed^n$ is always finite, and hence all the representing vectors can easily be cut out by a finite union of hyperplanes.
\end{rem}

\bigskip

\section{Notation and tools}
\label{tools}

We start with some notation and technical lemmas, following~\cite{cfh}. Let $K$ be a number field and $d := [K:\que]$ the global degree of $K$ over $\que$. Let $M(K)$ be the set of all places of $K$. For each place $v \in M(K)$, we write $K_v$ for the completion of $K$ at $v$, and let $d_v := [K_v : \que_v]$ be the local degree of $K$ at $v$. For each place $v \in M(K)$ we define the absolute value $|\ |_v$ to be the unique absolute value on $K_v$ that extends either the usual absolute value on $\real$ or $\cee$ if $v | \infty$, or the usual $p$-adic absolute value on $\que_p$ if $v|p$, where $p$ is a prime. For each non-zero $a \in K$ the following product formula is satisfied:
\begin{equation}
\label{product_formula}
\prod_{v \in M(K)} |a|^{d_v}_v = 1, \quad \mbox{ for all $a \in K^\times$}.
\end{equation}
For each $v \in M(K)$, we define a local height $H_v$ on $K_v^n$ by
$$H_v(\bx) = \max_{1 \leq i \leq n} |x_i|^{d_v}_v,\quad \mbox{ for each $\bx \in K_v^n$},$$
and take a product to define a global height function on $K^n$:
\begin{equation}
\label{global_heights}
H(\bx) = \left( \prod_{v \in M(K)} H_v(\bx) \right)^{1/d}
\end{equation}
for each $\bx \in K^n$. This height function is {\it homogeneous}, in the sense that it is defined on the projective space over $K^n$, thanks to the product formula \eqref{product_formula}. We also define the {\it inhomogeneous} height
$$h(\bx) = H(1,\bx),$$
which generalizes the Weil height on algebraic numbers. Clearly, $h(\bx) \geq H(\bx)$ for each $\bx \in K^n$. We extend the height functions $H$ and $h$ to polynomials by evaluating the height of their coefficient vectors. Finally, we define the height of an $m$-dimensional subspace $V \subseteq K^n$ as
$$H(V) := H(\bx_1 \wedge \dots \wedge \bx_m),$$
where $\bx_1,\dots,\bx_m$ is a basis for $V$ and $\bx_1 \wedge \dots \wedge \bx_m$ is viewed as a vector in $K^{\binom{n}{m}}$ under the standard Pl\"ucker coordinate embedding. Due to the product formula, this global height does not depend on the choice of the basis, and $H(K^n)=1$.
\smallskip

Now we review some basic tools that we will need for our main arguments.

\begin{lem} [\cite{cfh}, Lemma~2.1] \label{sum_height} For $\xi_1,...,\xi_\ell \in K$ and $\bx_1,...,\bx_\ell \in K^n$,
$$H \left( \sum_{i=1}^\ell \xi_i \bx_i \right) \leq h \left( \sum_{i=1}^\ell \xi_i \bx_i \right) \leq \ell h(\bxi) \prod_{i=1}^\ell h(\bx_i),$$
where $\bxi = (\xi_1,...,\xi_\ell) \in K^\ell$.
\end{lem}

\begin{lem} \label{H_Z} Let $\bx = (x_1, \ldots, x_n) \in \zed^n$, then $h(\bx) = \vert \bx \vert$, where $|\bx|$ is the sup-norm of $\bx$, and $H(\bx) = h(\bx)$ if $\bx$ is primitive, i.e. $\gcd(x_1, \ldots, x_n) = 1$.  Moreover, for $\bx,\bwy \in \zed^n$, $h(\bx + \bwy) \leq h(\bx) + h(\bwy)$.
\end{lem}

\proof
Let $\bx \in \zed^n$.  Then $H_p(\bx, 1) = 1$ for all finite primes $p$ and $H_\infty(\bx, 1) = \max\{\vert x_1\vert, \ldots, \vert x_n\vert \} $ which is the sup-norm $\vert \bx \vert$.  Therefore, $h(\bx) = H(\bx, 1) = \vert x \vert$.  Furthermore, $|\bx| = H_{\infty}(\bx)$, while the product over $p$-adic absolute values
$$\prod_{v \nmid \infty} H_v(\bx) = \gcd(x_1,\dots,x_n)^{-1}.$$
Thus, $H(\bx) = \vert \bx \vert/\gcd(x_1, \ldots, x_n)$ and hence $H(\bx) = h(\bx)$ if $\bx$ is primitive.  The last assertion is clear because $h$ is just the sup-norm.
\endproof

\begin{lem} [\cite{cfh}, Lemma~2.2] \label{lem_4.7} Let $V$ be a subspace of $K^n$, $n \geq 2$, and let subspaces $U_1,\dots,U_k \subseteq V$ and vectors $\bx_1,\dots,\bx_{\ell} \in V$ be such that
$$V = \spn_K \{ U_1,\dots,U_k,\bx_1,\dots,\bx_{\ell} \}.$$
Then
$$H(V) \ll H(U_1) \cdots H(U_k) H(\bx_1) \cdots H(\bx_{\ell}),$$
where the constant in the upper bound depends on $n$ and $\ell$.
\end{lem}

\begin{lem} [\cite{cfh}, Lemma 2.4] \label{intersect} Let $U_1, U_2$ be subspaces of $K^n$, then
$$H(U_1 \cap U_2) \leq H(U_1) H(U_2).$$
\end{lem}

\begin{lem} \label{dual} Let $Q$ be a nonzero quadratic form in $n$ variables over $K$ and let $V \subsetneq K^n$ be an $m$-dimensional subspace, $1 \leq m < n$.  Then
$$H(\perp_Q(V)) \ll H(Q)^m H(V),$$
where the constant in the upper bound depends on $n$, $m$, and $K$.
\end{lem}

\proof
Let $\bx_1,\dots,\bx_m$ be a small-height basis for $V$ guaranteed by Siegel's lemma (see~\cite{vaaler:siegel}), so
\begin{equation}
\label{siegel}
\prod_{i=1}^m h(\bx_i) \ll H(V),
\end{equation}
where the constant in the upper bound depends on $K$ and $m$. By Lemma~2.3 of~\cite{cfh},
\begin{equation}
\label{2.3}
H(\perp_Q(V)) \ll H(Q)^m \prod_{i=1}^m H(\bx_i),
\end{equation}
where the constant in the upper bound depends only on $n$ and $m$. The conclusion follows by combining~\eqref{siegel} with~\eqref{2.3}.
\endproof

We conclude this section by presenting the proof of Corollary \ref{smallrepV}.  From now on, for a string of parameters $I$,  we will write ``$f \ll_I g$" to mean that $f \leq C\, g$ for some constant $C$ depending only on the parameters in $I$.
\bigskip

\proof[Proof of Corollary \ref{smallrepV}]
As in the proof of Lemma \ref{dual}, $V$ has a basis $\bwy_1,\dots,\bwy_m$ satisfying
\begin{equation} \label{siegelZ}
\prod_{i=1}^m h(\bwy_i) \ll_m H(V).
\end{equation}
But the version of Siegel's lemma in \cite[Theorem 9]{vaaler:siegel} says more: $\bwy_1, \ldots, \bwy_m$ can be chosen to be in $V \cap \mathbb Z^n$.  For every $\bwy \in V\cap \mathbb Z^n$, $h(\bwy)$ is simply the sup-norm of the integral vector $\bwy$ which is comparable to the Euclidean length of $\bwy$.  Then by the Hadamard Inequality~\cite[Theorem 2.1.1]{martinet}, the Hermite Inequality~\cite[Theorem 2.2.1]{martinet}, and \eqref{siegelZ} together, the $\mathbb Z$-lattice $V \cap \mathbb Z^n$ has a basis $\bx_1, \ldots, \bx_m$ such that
\begin{equation} \label{basisVZ}
\prod_{i=1}^m h(\bx_i) \ll_m H(V).
\end{equation}

Let $T$ be the $n \times m$ matrix $[\bx_1\, \cdots \, \bx_m]$ and $Q'(X_1, \ldots, X_m)$ be the integral quadratic form with Gram matrix $T^t Q T$.  Then $Q'$ is a regular integral quadratic form in $m$ variables and $t$ is integrally represented by $Q'$.  By Theorem \ref{smallrep}, there exists $\bu \in \mathbb Z^m$ such that $Q'(\bu) = t$ and $h(\bu) \ll_m h(Q'_t)^{\ell(m)}$.  If $\bu = (u_1, \ldots, u_m)$, then $\bz: = u_1\bx_1 + \cdots + u_m\bx_m$ is a vector in $V\cap \mathbb Z^n$ such that $Q(\bz) = t$.

By Lemma \ref{sum_height} and \eqref{basisVZ}, we have
\begin{equation} \label{u-height}
h(\bz) \ll_m h(\bu) \prod_{i=1}^m h(\bx_i) \ll_m h(\bu) H(V) \ll_m h(Q'_t)^{\ell(m)} H(V).
\end{equation}
Since the matrices $T$ and $Q$ have integral entries, \eqref{basisVZ} implies that
\begin{equation} \label{Q-height}
h(Q'_t) \ll_m h(Q_t)h(T)^2 \leq h(Q_t) \left( \prod_{i=1}^m h(\bx_i)\right)^2 \ll_m h(Q_t)H(V)^2,
\end{equation}
where by $h(T)$ we mean height of $T$ viewed as a vector. Combining \eqref{u-height} and \eqref{Q-height} we obtain \eqref{smallrepVbound} in the Corollary.
\endproof

\bigskip

\section{Proof of Theorem~\ref{QZ}}
\label{field}

Starting with an algebraic set $\Z_{\bs} \subset K^n$ as in~\eqref{Z_K}, let us consider its homogenization $\Z_{\bs^*} \subset K^{n+1}$ as in~\eqref{Z*_K}. Given the $m$-dimensional subspace $V$ of $K^n$, let $V^*=V \oplus K\be_{n+1} \subseteq K^{n+1}$, where $\be_{n+1}$ is the last standard basis vector in~$K^{n+1}$. Then $V^*$ is an $(m+1)$-dimensional subspace of $K^{n+1}$, and by Lemma~\ref{lem_4.7},
\begin{equation}
\label{HV*}
H(V^*) \ll_n H(V) H(\be_{n+1}) = H(V),
\end{equation}
since the height of $V$ does not change after the embedding $V \to (V,0)$ into $K^{n+1}$. Let $t \neq 0$ be an element of $Q(V)$, and $Q_t^* \in K[X_1, \ldots, X_{n+1}]$ be the homogenization of the quadratic polynomial $Q(X_1, \ldots, X_{n+1}) - t$, i.e.
$$Q_t^*(X_1,\dots,X_{n+1}) = Q(X_1,\dots,X_n) - t X_{n+1}^2.$$
Then for any $\bx \in V$, $Q(\bx) = t$ if and only if $Q_t^*(\bx,1) = 0$. The hypothesis of the theorem implies that there exists $\bz \in V \setminus \Z_{\bs}$ such that $Q(\bz) = t$. Let
$$\Z'_{\bs} = \Z_{\bs^*} \cup \{ (\bx, 0) : \bx \in K^n \} \subset K^{n+1}.$$
Then $Q_t^*(\bz,1)=0$ and $(\bz,1) \in V^* \setminus \Z'_{\bs}$, meaning that the quadratic form $Q_t^*$ on $K^{n+1}$ has a nontrivial zero in $V^* \setminus \Z'_{\bs}$. Notice also that for every maximal totally isotropic subspace $W$ of $(V,Q)$, $(W,0)$ is a totally isotropic subspace of $(V^*,Q_t^*)$, and hence the dimension $w^*$ of a maximal totally isotropic subspace of $(V^*,Q_t^*)$ is at least $w$. Further, if $w=0$, i.e. if $(V,Q)$ is anisotropic, then $w^*=1$. Then, by Theorem~\ref{miss_hyper}, there must exist $\bwy \in V^* \setminus \Z'_{\bs}$ such that $Q_t^*(\bwy) = 0$ and
$$H(\bwy) \ll_{m, M_{\bs}} H(Q_t^*)^{\frac{(m+1)-w^*+1}{2}} H(V^*)^2 \ll_n H(Q_t^*)^{\frac{m-w+2}{2}} H(V)^2,$$
by~\eqref{HV*}. Since $\bwy \notin \Z'_{\bs}$, $y_{n+1} \neq 0$, so define
$$\bwy' = \frac{1}{y_{n+1}} \bwy = \left( \frac{y_1}{y_{n+1}},\dots,\frac{y_n}{y_{n+1}},1 \right),$$
and let $\bz = \left( \frac{y_1}{y_{n+1}},\dots,\frac{y_n}{y_{n+1}} \right) \in V$, i.e. $\bwy' = (\bz,1)$. Then, by the product formula,
\begin{equation}
\label{ht_z}
h(\bz) = H(\bwy') = H(\bwy) \ll_{m,n, M_{\bs}} H(Q_t^*)^{\frac{m-w+2}{2}} H(V)^2,
\end{equation}
and $Q_t^*(\bwy') = Q_t^*(\bwy) = 0$, which means that $Q(\bz) = t$.  This completes the proof of Theorem~\ref{QZ}.
\medskip

\begin{rem} Notice that the argument in the proof of Theorem~\ref{QZ} does not necessarily produce a point with coordinates in~$\O_K$, since we divided by~$y_{n+1}$. In the next section we produce weaker versions of such a result over~$\O_K$ when $K=\que$.
\end{rem}
\bigskip

\section{Proof of Theorem~\ref{subspaces}}
\label{proofs}

In this section we give the proof of Theorem~\ref{subspaces}.  We start with the following lemma, which is essentially Case 1 in the proof of the main theorem in \cite{dietmann}.  We extract its proof from \cite{dietmann} and present it here for the convenience of the reader.

\begin{lem} \label{case1}
Let $K$ be a number field and $Q(X_1, \ldots, X_n) \in K[X_1, \ldots, X_n]$ be a nonzero quadratic form.  Suppose that $(V, Q)$ is a regular isotropic subspace of $K^n$ of dimension $\geq 3$.  If $W_1, \ldots, W_k$ are hyperplanes of $V$, then $V$ has an isotropic vector outside of the union $\bigcup_{i=1}^k W_i$.
\end{lem}

\proof After fixing a choice of basis for $V$, we may assume that $V = K^n$ with $n \geq 3$.  For $i = 1, \ldots, k$, let $L_i$ be a nonzero linear form in $K[X_1, \ldots, X_n]$ whose kernel is $W_i$.  Let $\bu$ be an isotropic vector in $K^n$.  For every $\bx \in K^n$, let
$$T_\bu(\bx) = Q(\bx) \bu - 2 Q(\bx, \bu)\bx.$$
It is straightforward to check that $Q(T_\bu(\bx)) = 0$ for all $\bx \in K^n$.

Fix an index $i$.  We claim that $L_i(T_\bu(X_1, \ldots X_n))$ is not the zero polynomial.  Assume the contrary that $L_i(T_\bu(X_1, \ldots, X_n))$ is indeed the zero polynomial.  Suppose first that $L_i(\bu) \neq 0$.  Then
$$Q(\bx) = \frac{2Q(\bu, \bx)L_i(\bx)}{L_i(\bu)}$$
for all $\bx \in K^n$, which means that $Q(X_1, \ldots, X_n)$ is a product of two linear forms which is impossible since $(K^n, Q)$ is a regular quadratic space of dimension $\geq 3$.  Now, if $L_i(\bu) = 0$, then $Q(\bu, \bx)L_i(\bx) = 0$ for all $\bx \in K^n$, which is impossible again since both $L_i$ and $\bx \longmapsto Q(\bu, \bx)$ are not the zero linear functional on $K^n$.

So, the polynomial
$$F(X_1, \ldots, X_n): = L_1(T_\bu(X_1, \ldots, X_n)) \cdots L_k(T_\bu(X_1, \ldots, X_n))$$
is not the zero polynomial in $K[X_1, \ldots, X_n]$.  We can then take a vector $\bx \in K^n$ such that $F(\bx) \neq 0$.  Then $T_\bu(\bx)$ is an isotropic vector in $K^n$ outside of $\bigcup_{i=1}^k W_i$.
\endproof

\begin{lem} \label{Qy} Let $Q(X_1, \ldots, X_n) \in \zed[X_1, \ldots, X_n]$ be an integral quadratic form. Let $V \subseteq \que^n$ be an $m$-dimensional subspace, $3 \leq m \leq n$, such that the quadratic space $(V,Q)$ is regular and isotropic of Witt index $w$.  Let $W_1, \ldots, W_k$ be distinct hyperplanes of $V$. Then for every anisotropic vector $\bwy \in V\cap \zed^n$, there exists a point $\bz \in \left( V \cap \zed^n \right) \setminus \bigcup_{i = 1}^k W_i$ such that $Q(\bz) = Q(\bwy)$ and
$$h(\bz) \ll_{m,n,k} \begin{cases}
 h(\bwy)^{1 + \frac{2}{m-2}}\, h(Q)^{m + 2 + \frac{m + 4}{m-2}}\, H(V)^{4 + \frac{6}{m-2}} & \mbox{ if $w = 1$};\\
h(\bwy)^2\, h(Q)^{\frac{m-w+5}{2}}\, H(V)^2 & \mbox{ if $w \geq 2$}.
\end{cases}$$
\end{lem}

\proof
If $\bwy \not \in \bigcup_{i=1}^k W_i$, then we are done; so we assume from now on that $\bwy \in \bigcup_{i=1}^k W_i$.
\bigskip

\noindent{\em Case 1}:  Assume  that $w$, the Witt index of $V$, is at least 2.   Let $V_{\bwy}$ be the orthogonal complement of $\bwy$ in $V$, which is an $(m-1)$-dimensional regular isotropic subspace of $V$ such that $V = \mathbb Q[\bwy] \perp V_{\bwy}$.  We have two cases to consider, according to whether $V_{\bwy}$ is one of the $W_i$ or not.

Suppose first that $V_{\bwy}$ is not equal to any of the $W_i$.  Then $V_{\bwy} \cap W_i$ is a hyperplane of $V_{\bwy}$ for each $i$.  By Lemma \ref{case1}, $V_{\bwy}\setminus \bigcup_{i=1}^k W_i$ has an isotropic vector.  Then, by Theorem \ref{miss_hyper} and Lemmas \ref{lem_4.7}, \ref{intersect}, \ref{dual} there exists an isotropic vector $\bu \in (V_{\bwy}\cap \zed^n)\setminus \bigcup_{i=1}^k W_i$ with
\begin{equation} \label{bound-uy}
h(\bu) \ll_{m,n,k} h(Q)^{\frac{m-w+1}{2}}H(V_{\bwy})^2 \ll_{m,n,k} h(\bwy)^2 h(Q)^{\frac{m-w+5}{2}} H(V)^2.
\end{equation}
Note that $H(Q) = h(Q)$ because $Q$ has integer coefficients.  For $1 \leq r \leq k+1$,   let $\bz_r: = \bwy + r\bu$ which is in $(V\cap \zed^n)$.  Suppose that all of these $\bz_r$ are in $\bigcup_{i=1}^k W_i$.  By the Pigeon-Hole Principle, there must be $1 \leq r \neq s \leq k+1$ and $1 \leq j \leq k$  such that $\bz_r$ and $\bz_s$ are in $W_j$.  Then $\bz_r - \bz_s$ is in $W_j$, implying that $\bu$ is also in $W_j$ which is a contradiction.  Thus, there must be an $r$ such that $\bz: = \bz_r$ is not in $\bigcup_{i=1}^k W_i$ and $Q(\bz) = Q(\bwy)$.   Moreover,
\begin{equation} \label{bound-zy}
h(\bz)\leq h(\bwy) + r h(\bu) \ll_{m,n, k} h(\bwy)^2 h(Q)^{\frac{m-w+5}{2}} H(V)^2.
\end{equation}

Now, suppose that $V_{\bwy}$ is one of the $W_i$, say $V_{\bwy} = W_1$.  Then $V_{\bwy} \cap W_i$ is a hyperplane in $V_{\bwy}$ for every $i \geq 2$.  As is done in the last paragraph, we have an isotropic vector $\bu \in (V_{\bwy}\cap \zed^n)\setminus \bigcup_{i=2}^k W_i$ and $h(\bu)$ is bounded above as in \eqref{bound-uy}.   Since $\bu$ is in $W_1$ but $\bwy$ is not, all of the $\bz_r$ defined as in the last paragraph are not in $W_1$.  Suppose that all of these $\bz_r$ are in $\bigcup_{i=2}^k W_i$.   As is argued before, there will be two different indices $r, s \in \{1, \ldots k + 1\}$ and $j \in \{2, \ldots, k\}$ such that $\bz_r - \bz_s \in W_j$.  This means that $\bu$ is in $W_j$ for some $j \geq 2$, which is impossible.  We may now argue as in the last paragraph to obtain an $\bz \in (V\cap \zed^n)\setminus \bigcup_{i=1}^k W_i$ with $Q(\bz) = Q(\bwy)$ and $h(\bz)$ bounded above as in \eqref{bound-zy}.
\bigskip

\noindent {\em Case 2}: Suppose that $w = 1$.  The space $V$ itself is regular and isotropic.  By Lemma \ref{case1} and Theorem \ref{miss_hyper}, there exists an isotropic vector $\bu \in (V\cap \zed^n)\setminus \bigcup_{i=1}^k W_i$ with
\begin{equation} \label{bound-u}
h(\bu) \ll_{m,n,k} h(Q)^{\frac{m}{2}}H(V)^2.
\end{equation}
Suppose first that $Q(\bu, \bwy) = 0$.   We could argue as in Case 1 and obtain a vector $\bz \in (V\cap \zed^n)\setminus \bigcup_{i=1}^k W_i$ such that $Q(\bz) = Q(\bwy)$ and
\begin{equation} \label{bound-z11}
h(\bz)\leq h(\bwy) + r h(\bu) \ll_{m,n, k} h(\bwy) h(Q)^{\frac{m}{2}} H(V)^2.
\end{equation}
Note that this upper bound is better than the one in the statement of the lemma.

Now, let us suppose that $Q(\bu, \bwy) \neq 0$.   Then $\mathbb Q[\bu, \bwy]$ is a 2-dimensional regular subspace of $V$, and its orthogonal complement $G$ in $V$ is an $(m-2)$-dimensional anisotropic subspace of $V$.  By Siegel's lemma \cite[Theorem 9]{vaaler:siegel}, there  exist vectors $\bx_1, \ldots, \bx_s$ in $\zed^n$ which form a basis for $G$ and
\begin{equation} \label{bound-xproduct}
\prod_{i=1}^{m-2} h(\bx_i) \ll_m H(G).
\end{equation}
Let $\bx$ be a vector among the $\bx_i$ that has the smallest height.  Then, combining this with \eqref{bound-u} and using Lemmas \ref{lem_4.7}, \ref{intersect}, and~\ref{dual}, we have
\begin{equation}\label{bound-x}
h(\bx) \ll_m H(G)^{\frac{1}{m-2}}\ll_{m,n} h(\bwy)^{\frac{1}{m-2}}\, h(Q)^{\frac{m+4}{2(m-2)}}\, H(V)^{\frac{3}{m-2}}.
\end{equation}
This $\bx$ is necessarily an anisotropic vector in $V\cap \zed^n$.  By replacing $\bx$ with $2\bx$ if necessary, we may assume that $Q(\bx)$ is a nonzero even integer.

For each integer $\ell \geq 1$, let $a_\ell = 2^{\ell - 1} + \cdots + 2 + 1$ and
$$I_\ell = \{(k + 1)a_\ell - k,  \ldots, (k + 1)a_\ell\}.$$
For each $r \in I_\ell$, let
$$\bz_r: = \bwy + rQ(\bwy, \bu)\bx - \frac{r^2}{2}Q(\bx)Q(\bwy, \bu)\bu.$$
It is easy to check that $\bz_r \in V\cap \zed^n$ and  $Q(\bz_r) = Q(\bwy)$.  Suppose that for all $\ell \in \{1, \ldots, k+1\}$ and all $r \in I_\ell$, $\bz_r$ are in $\bigcup_{i=1}^k W_i$.  Fix an $\ell$; by the Pigeon-Hole Principle, there must be two distinct indices $r, s \in I_\ell$ and $i(\ell) \in \{1, \ldots, k\}$ such that both $\bz_r$ and $\bz_s$ are in $W_{i(\ell)}$.  Then $\bz_r - \bz_s$ is in $W_{i(\ell)}$ and hence
$$\bx - \frac{1}{2}Q(\bx)(r + s)\bu \in W_{i(\ell)}.$$
By virtue of the Pigeon-Hole Principle one more time,  there must be two different indices $\ell_1, \ell_2 \in \{1, \ldots, k + 1\}$, $i_0 \in \{1, \ldots, k\}$, and pairs of distinct integers $r_1, s_1 \in I_{\ell_1}$ and $r_2, s_2 \in I_{\ell_2}$ such that both $\bz_{r_1} - \bz_{s_1}$ and $\bz_{r_2} - \bz_{s_2}$ are in $W_{i_0}$.  Then,
$$\frac{1}{2}Q(\bx)\left( (r_2 + s_2) - (r_1 + s_1)\right)\bu \in W_{i_0}.$$
Without loss of generality, we may assume that $\ell_1 < \ell_2$.  Then
\begin{eqnarray*}
r_1 + s_1 & \leq & 2(k+1)a_{\ell_1}\\
    & = & (k+1)a_{\ell_1 + 1} - (k+1)\\
    & \leq & (k+1)a_{\ell_2} - k  - 1\\
    & <  & 2(k+1)a_{\ell_2} - 2k\\
    & \leq & r_2 + s_2,
\end{eqnarray*}
and hence $\bu \in W_{i_0}$, which is impossible. Thus, there must be one of these $\bz_r$  in $(V\cap \zed^n)\setminus \bigcup_{i=1}^k W_i$.  Let $\bz$ be this $\bz_r$.   A crude estimate shows that
\begin{equation} \label{bound-z2}
h(\bz) \ll_{m,n,k} h(\bwy)h(Q)^2\, h(\bu)^2\, h(\bx)^2.
\end{equation}
Using \eqref{bound-u}, \eqref{bound-x}, and \eqref{bound-z2} we obtain
\begin{equation} \label{bound-z3}
h(\bz) \ll_{m,n,k} h(\bwy)h(\bwy)^{\frac{2}{m-2}}\, h(Q)^{m+2 + \frac{m+4}{m-2}}\, H(V)^{4 + \frac{6}{m-2}}.
\end{equation}
The lemma is proved by combining \eqref{bound-zy}, \eqref{bound-z11}, and \eqref{bound-z3}.
\endproof

\proof[Proof of Theorem~\ref{subspaces}] The theorem now follows upon combining Corollary~\ref{smallrepV} and Lemma~\ref{Qy}, noticing that $h(Q) \leq h(Q_t)$.
\endproof
\bigskip

{\bf Acknowledgement:} We thank the referee for a very careful and thorough reading and many helpful comments which have improved the quality of the paper.
\bigskip

\bibliographystyle{plain}  
\bibliography{small_representations-1}        

\end{document}